\documentclass[english]{amsart}

\usepackage{esint}
\usepackage[svgnames]{xcolor} 
\usepackage[colorlinks,citecolor=red,pagebackref,
hypertexnames=false,breaklinks]{hyperref}
\usepackage{pgf,tikz}
\usepackage{pdfsync}

\usepackage{dsfont}
\usepackage{url}
\usepackage[utf8]{inputenc}
\usepackage[T1]{fontenc}
\usepackage{lmodern}
\usepackage{babel}
\usepackage{mathtools}  
\usepackage{amssymb}
\usepackage{lipsum}
\usepackage{mathrsfs}
\usepackage{color}

\usepackage[skip=2.2pt plus 1pt, indent=12pt]{parskip}

\newtheorem{theorem}{Theorem}[section]

\newtheorem{corollary}{Corollary}[section]

\numberwithin{equation}{section}

\title[Stable determination of initial data]{Stable determination of the initial data in an IBVP for the wave equation outside a non-trapping obstacle}

\author[Mourad Choulli]{Mourad Choulli}
\address{Universit\'e de Lorraine}
\email{mourad.choulli@univ-lorraine.fr}

\begin{document}

\begin{abstract}
We establish a double logarithmic stability inequality for the problem of determining the initial data in an IBVP for the wave equation outside a non-trapping obstacle from two localized measurements.
\end{abstract}



\maketitle

%\tableofcontents

\section{Introduction}\label{section1}

Let $n\ge 3$ be an odd integer and let $K$ be the closure of a star-shaped bounded open subset of $\mathbb{R}^n$ with $C^\infty$ boundary. Let $\Omega =\mathbb{R}^n\setminus K$. Throughout this text we use the notation $\mathbb{R}_+=]0,\infty[$.

Define the unbounded operator $A: L^2(\Omega)\rightarrow L^2(\Omega)$ as follows
\[
Au=\Delta u\quad \mbox{and}\quad D(A)=H_0^1(\Omega)\cap H^2(\Omega).
\]
It is a classical result that $-A$ is self-adjoint positive operator.

We endow the Hilbert space $\mathcal{H}=H_0^1(\Omega)\times L^2(\Omega)$ with the scalar product
\[
\langle (u_1,v_1)|(u_2,v_2)\rangle =\int_\Omega [\nabla u_1\cdot \nabla \overline{u_2}+v_1\overline{v_2}]dx,\quad (u_j,v_j)\in \mathcal{H},\; j=1,2.
\] 
Let $\mathcal{A}:\mathcal{H}\rightarrow \mathcal{H}$ be the unbounded operator given as follows
\[
\mathcal{A}(u,v)=(v,Au)\quad \mbox{and}\; D(\mathcal{A})=(H_0^1(\Omega)\cap H^2(\Omega))\times H_0^1(\Omega).
\]
As $-A$ is self-adjoint positive operator, $\mathcal{A}$ is m-dissipative and therefore, according to Lumer-Phillips theorem, $\mathcal{A}$ generates a contraction semigroup (e.g. \cite[Theorem 3.8.4]{TW}). Whence, for every $(f,g)\in (H_0^1(\Omega)\cap H^2(\Omega))\times H_0^1(\Omega)$ the IBVP 
\begin{equation}\label{1.2}
\partial_t^2u-\Delta u=0\; \mbox{in}\; \mathbb{R}_+\times \Omega,\quad u_{|\mathbb{R}_+\times \partial \Omega}=0,\quad (u(0,\cdot),\partial_tu(0,t))=(f,g)
\end{equation}
has a unique solution $u=u(f,g)$ so that $u\in C(\overline{\mathbb{R}}_+,H_0^1(\Omega)\cap H^2(\Omega))$, $\partial_tu\in C(\overline{\mathbb{R}}_+,H_0^1(\Omega))$ and $\partial_t^2u\in C(\overline{\mathbb{R}}_+,L^2(\Omega))$. Furthermore, we have
\begin{align}
\|u(t,\cdot)\|_{H^2(\Omega)}+\|\partial_tu(t,\cdot)\|_{H^1(\Omega)}+&\|\partial_t^2u(t,\cdot)\|_{L^2(\Omega)}. \label{1.3}
\\
&\le \mathfrak{m}\left( \|f\|_{H^2(\Omega)}+\|g\|_{H^1(\Omega)}\right),\quad t\ge 0,\nonumber
\end{align}
where $\mathfrak{m}=\mathfrak{m}(\Omega)>0$ is a constant.

Let $a_j<b_j$, $1\le j\le n$, be chosen in such a way that $\displaystyle Q_0:=\prod_{j=1}^n[a_j,b_j]\subset \Omega$. Let $\Sigma \supset K$ be a bounded domain of $\mathbb{R}^n$ with $C^\infty$ boundary so that $Q=\Sigma \setminus K\Supset Q_0 $.  In addition, we assume that there exists $Q_1\subset Q$ satisfying $Q_0\Subset Q_1$ and $Q\setminus Q_1$ is connected. Set
\[
\mathbf{D}=\{ (f,g)\in (H_0^1(\Omega)\cap H^2(\Omega))\times H_0^1(\Omega);\; \mbox{supp}(f)\cup \mbox{supp}(g)\subset Q_0\}.
\]

In light of the results established in \cite[Section 3]{MRS} we can state the following theorem.
\begin{theorem}\label{theorem1}
There exist $\varkappa=\varkappa (\Omega,Q) >0$ and $\delta=\delta (\Omega,Q) >0$ so that for any $(f,g)\in \mathbf{D}$ we have
\begin{align}
&\|u(t,\cdot)\|_{H^2(Q)}+\|\partial_tu(t,\cdot)\|_{H^1(Q)}+\|\partial_t^2u(t,\cdot)\|_{L^2(Q)}\label{1.4}
\\
&\hskip 5cm \le \varkappa e^{-\delta t}\left(\|f\|_{H^2(\Omega)}+\|g\|_{H^1(\Omega)}\right),\quad t\ge 0,\nonumber
\end{align}
where $u=u(f,g)$.
\end{theorem}

It should be observed that $u=u(f,g)$ given by theorem \ref{theorem1.1} satisfies $u\in L^1((\mathbb{R}_+,t^kdt), H^ 2(\Omega))$ and
\begin{align*}
&\left(\int_0^{+\infty}u_{|Q}t^kdtdx\right)_{\big|\partial Q}=\int_0^{+\infty}u_{|\partial Q}t^kdtdx,
\\
&\nabla  \left(\int_0^{+\infty}u_{|Q}t^kdtdx\right)_{\big|\partial Q}=\int_0^{+\infty}\nabla u_{|\partial Q}t^kdtdx.
\end{align*}
The notation $|\partial Q$ is to be understood in the usual sense of trace. 

In what follows, $S$ denotes an arbitrary fixed nonempty open subset of $\partial \Sigma$. Let $(f,g)\in \mathbf{D}$ and $u=u(f,g)$. The identities above allow us to define 
\[
N(u)=\sum_{j=0,1}\left(\|u\|_{L^1((\mathbb{R}_+,t^jdt),L^2(S))}+\|\nabla u\|_{L^1((\mathbb{R}_+,t^jdt),L^2(S))}\right).
\]

Let $\mathbf{D}_\vartheta=\{(h,k)\in \mathbf{D};\; \|h\|_{H^2(\Omega)}+\|k\|_{H^1(\Omega)}\le \vartheta\}$, $\vartheta >0$.

\begin{theorem}\label{theorem1.0}
Let $\vartheta >0$. There exist $C=C(n,Q,S)$ and $\mathfrak{p}=\mathfrak{p}(n,Q,Q_0)$ satisfying : for any $(f,g),(\tilde{f},\tilde{g})\in \mathbf{D}_\vartheta$ so that $f-\tilde{f}$ and $g-\tilde{g}$ are independent of one of the variables $x_1,\ldots ,x_n$  we have
\begin{equation}\label{mi}
C\left(\|f-\tilde{f}\|_{H^1(\Omega)}+\|g-\tilde{g}\|_{L^2(\Omega)}\right)\le \rho^{-1}\vartheta +e^{e^{\mathfrak{p}\rho}}N(u-\tilde{u}),\quad \rho\ge 1,
\end{equation}
where $u=u(f,g)$ and $\tilde{u}=u(\tilde{f},\tilde{g})$.
\end{theorem}

Note that $N(u-\tilde{u})$ in Theorem \ref{theorem1.0} can be replaced by 
\[
\tilde{N}(u-\tilde{u})=\sum_{j=0,1}\left(\|u-\tilde{u}\|_{L^1((\mathbb{R}_+,t^jdt),H^1(S))}+\|\partial_\nu (u-\tilde{u})\|_{L^1((\mathbb{R}_+,t^jdt),L^2(S))}\right),
\]
where $\nu$ denotes the unitary exterior normal vector field at $\partial Q$.

For clarity, we intentionally choose to state our main result with non-optimal assumptions: $K$ is star-shaped and the speed of sound in the wave equation \eqref{1.2} is assumed identically equal to $1$. More general situations will be discussed and commented in Section \ref{section3}. We also provide in Section \ref{section3} some references related to the problem we consider in this article. 

Our study is motivated by thermo-acoustic and photo-acoustic tomography.

The proof of  Theorem \ref{theorem1.0} relies on the exponential decay of the local energy of the wave equation outside a non-trapping obstacle in odd space dimensions (Theorem \ref{theorem1.1})  combined with a quantized version of the unique continuation from an arbitrary Cauchy datum for elliptic equations.

\section{Proof of Theorem \ref{theorem1.0}}\label{section2}

In what follows we use the notation $\Pi_\mu =\{z\in \mathbb{C};\; \Re z >-\mu\}$, $\mu >0$. Let $(f,g)\in \mathbf{D}$ and $u=u(f,g)$. Referring to Theorem \ref{theorem1}, we can define
\begin{equation}\label{1.5}
v(z,\cdot )=\int_0^{+\infty}e^{-tz}u(t,\cdot)_{|Q}dt,\quad z\in \Pi_\delta.
\end{equation}
where $\delta=\delta(\Omega,Q)$ is the same as in Theorem \ref{theorem1}. Note that $z\mapsto v(z,\cdot )$ is nothing other than the Laplace transform of $t\mapsto u(t,\cdot)_{|Q}$.
Using \eqref{1.4}, we show that $v\in \mathcal{H}(\Pi_\delta, H^2(Q))$, where $\mathcal{H}(\Pi_\delta, H^2(Q))$ is the space of holomorphic functions defined on $\Pi_\delta$ with values in $H^2(Q)$, and
\begin{equation}\label{1.6}
z^2v(z,x)-\Delta v(z,x)=zf(x)+g(x),\quad (z,x)\in \Pi_\delta \times Q.
\end{equation}
Observing that $\mathbb{D}_\delta =\{z\in \mathbb{C}; |z|<\delta\}\subset \Pi_\delta$, we write $v$ as an entire series in $z$:
\begin{equation}\label{1.1}
v(z,x)=\sum_{k\ge 0}v_k(x)z^k,\quad (z,x)\in \mathbb{D}_\delta\times Q,
\end{equation}
where $v_k\in H^2(Q)$, $k\ge 0$. Inserting this formula into \eqref{1.6}, we find
\begin{align}
&-\Delta v_0=g\quad \mbox{in}\; Q, \label{1.7}
\\
&-\Delta v_1=f\quad \mbox{in}\; Q, \label{1.8}
\\
&-\Delta v_k=-v_{k-2}\quad \mbox{in}\; Q,\; k\ge 2. \label{1.9}
\end{align}
Note that since $v_k$ is given by 
\[
v_k=\frac{1}{k!}\partial_z^kv(z,\cdot)_{|z=0},\quad k\ge 0,
\]
formula \eqref{1.5} yields
\begin{equation}\label{1.10}
v_k=\frac{(-1)^k}{k!}\int_0^{+\infty}u(t,\cdot)_{|Q}t^kdt,\quad k\ge 0.
\end{equation}

Let $\tilde{u}=u(\tilde{f},\tilde{g})$ with $(\tilde{f},\tilde{g})\in \mathbf{D}$. Define $\tilde{v}$ and $(\tilde{v}_k)$ by replacing in \eqref{1.5} and \eqref{1.1} $u$ with $\tilde{u}$.

Let $w_k=v_k-\tilde{v}_k$, $k\ge 0$. As $v_k{_{|\partial K}}=\tilde{v}_k{_{|\partial K}}=0$, $k\ge 0$, we obtain from inequalities \eqref{1.7}-\eqref{1.9} and inequalities \eqref{1.7}-\eqref{1.9} with $(v_k)$ replaced by $(\tilde{v}_k)$
\begin{align}
&-\Delta w_0=g-\tilde{g}\; \mbox{in}\;Q, \quad w_0{_{|\partial K}}=0, \label{1.14}
\\
&-\Delta w_1=f-\tilde{f}\; \mbox{in}\;Q, \quad w_1{_{|\partial K}}=0, \label{1.15}
\\
&-\Delta w_k=-w_{k-2}\; \mbox{in}\;Q, \quad w_k{_{|\partial K}}=0,\quad  k\ge 2. \label{1.16}
\end{align}

Let $\mathscr{H}(Q)=\{h\in H^2(Q);\; \Delta h=0\}$. Using Green's inequality, we get 
\begin{equation}\label{1.19}
\int_Q(g-\tilde{g})\phi dx=-\int_{\partial Q}\partial_\nu w_0\phi ds+\int_{\partial\Sigma}w_0\partial_\nu \phi ds,\quad \phi \in \mathscr{H}(Q).
\end{equation}

If $\Lambda =\{\xi=\eta +i\zeta;\; \eta,\zeta \in \mathbb{R}^n,\; |\eta|=|\zeta|,\; \eta \cdot \zeta =0\}$ we check that $e^{-ix\cdot \xi}\in \mathscr{H}(Q)$  for every $\xi \in \Lambda$. Let $\Lambda_0$ be the subset of $\Lambda$ consisting of $\xi=\eta +i\zeta$ given by
\[
\eta=(\eta',0)\in \mathbb{R}^{n-1} \times \mathbb{R},\quad \zeta =(0,\zeta_n)\in \mathbb{R}^{n-1} \times \mathbb{R},
\]
and  $|\eta'|=|\zeta_n|$.

Below we write $x=(x',x_n)\in \mathbb{R}^{n-1}\times \mathbb{R}$. Without loss of generality, we assume that $g-\tilde{g}=G(x')$.  Let $I_n$ the projection of $Q$ with respect to the variable $x_n$, $\displaystyle \gamma =\max_{I_n}|x|$, $\displaystyle \tau=\max_{[a_n,b_n]}|x_n|$, $\displaystyle P=\prod_{j=1}^{n-1}[a_j,b_j]$, $\eta '\in \mathbb{R}^{n-1}$ and $F$ the Fourier transform of $G\chi_P$. Then \eqref{1.19} with $\phi=e^{-ix\cdot \xi}$, $\xi \in \Lambda_0$, yields
\begin{equation}\label{1.20}
e^{-\tau |\eta'|}|F(\eta ')|\le e^{\gamma |\eta'|}\| \partial_\nu w_0\|_{L^1(\partial Q)}+\sqrt{2}|\eta'|e^{\gamma |\eta'|}\| w_0\|_{L^1(\partial \Sigma)}.
\end{equation}
Set
\[
N_0(u-\tilde{u})=\|\partial_\nu (u-\tilde{u})\|_{L^1(\mathbb{R}_+\times \partial Q)}+ \|u-\tilde{u}\|_{L^1(\mathbb{R}_+\times \partial \Sigma)}.
\]
Using \eqref{1.10}, we get from \eqref{1.20}
\[
|F(\eta ')|\le \sqrt{2}e^{\sigma |\eta '|}N_0(u-\tilde{u}),
\]
where $\sigma =1+\tau+\gamma$. Therefore we have  for any $\rho >0$
\begin{equation}\label{1.21}
\int_{\mathcal{B}_\rho}|F(\eta ')|^2d\eta '\le 2\rho^{n-1}|\mathcal{B}_1|e^{2\sigma \rho} [N_0(u-\tilde{u})]^2,
\end{equation}
where $\mathcal{B}_\rho$ is the ball of $\mathbb{R}^{n-1}$ with radius $\rho$ and center $0$. Thus we obtain
\begin{equation}\label{1.29}
\int_{\mathcal{B}_\rho}|F(\eta ')|^2d\eta '\le 2(n-1)!|\mathcal{B}_1|e^{(2\sigma+1) \rho} [N_0(u-\tilde{u})]^2.
\end{equation}
On the other hand, we have
\begin{equation}\label{1.22}
\int_{|\eta '|\ge \rho}|F(\eta ')|^2d\eta '\le \rho^{-2}\|g-\tilde{g}\|_{H^1(\mathbb{R}^{n-1})}^2.
\end{equation}
As $\|F\|_{L^2(\mathbb{R}^{n-1})}=\|g-\tilde{g}\|_{L^2(\mathbb{R}^{n-1})}$ by Parseval's inequality, \eqref{1.21} and \eqref{1.22} imply 
\[
\|g-\tilde{g}\|_{L^2(\mathbb{R}^{n-1})}\le \mathbf{c} \rho^{-1}\|g-\tilde{g}\|_{H^1(\mathbb{R}^{n-1})} +e^{\alpha \rho}N_0(u-\tilde{u}),
\]
where $\displaystyle \alpha=\sigma +\frac{1}{2}$ and $\mathbf{c}= \sqrt{2(n-1)!|\mathcal{B}_1|}$. Whence
\begin{equation}\label{1.27}
\|g-\tilde{g}\|_{L^2(\Omega)}\le \mathbf{c}  \rho^{-1}\|g-\tilde{g}\|_{H^1(\Omega)} +(b_n-a_n)^{\frac{1}{2}}e^{\alpha \rho}N_0(u-\tilde{u}).
\end{equation}
 
Now,  from \cite[Theorem 1.1]{Ch} applied to $w_0{_{|Q\setminus \overline{Q}_1}}$ with $\displaystyle s=\eta_0=\frac{1}{4}$, $\displaystyle \eta_1=\frac{1}{2}$ and the comments in \cite[Subsection 1.3]{Ch}, it follows that there exist $C=C(n,Q,S)>0$, $\mathfrak{b}=\mathfrak{b}(n,Q)>0$ and $\gamma=\gamma(n,Q)>0$ so that for every $0<r<1$ 
\[
CN_0(u-\tilde{u})\le r^{\frac{1}{8}}\|u-\tilde{u}\|_{L^1(\mathbb{R}_+,H^2(Q))}+e^{\mathfrak{b}r^{-\gamma}}N(u-\tilde{u}),
\]
which in combination with \eqref{1.4} yields
\begin{equation}\label{1.28}
CN_0(u-\tilde{u})\le r^{\frac{1}{8}}\left(\|f-\tilde{f}\|_{H^2(\Omega)}+\|g-\tilde{g}\|_{H^1(\Omega)}\right)+e^{\mathfrak{b}r^{-\gamma}}N(u-\tilde{u}).
\end{equation}

Let $\vartheta >0$ and assume that $(f,g),(\tilde{f},\tilde{g})\in \mathbf{D}_\vartheta$. Putting together \eqref{1.27} and \eqref{1.28}, we obtain for any $0<r<1$ and $\rho >0$
\begin{equation}\label{1.30}
C\|g-\tilde{g}\|_{L^2(\Omega)}\le \vartheta \rho^{-1}+e^{\alpha \rho}\left(\vartheta r^{\frac{1}{8}}+e^{\mathfrak{b}r^{-\gamma}}N(u-\tilde{u})\right).
\end{equation}
Let $\rho \ge 1$. Choosing in \eqref{1.30} $r$ in such a way that $r^{-1/8}= \rho e^{\alpha \rho}$, we get
\[
C\|g-\tilde{g}\|_{L^2(\Omega)}\le \vartheta \rho^{-1}+e^{e^{\mathfrak{p}\rho}}N(u-\tilde{u}),
\]
where $\mathfrak{p}=\mathfrak{p}(n,Q,Q_0)>0$ is a constant.

To complete the proof, we use \eqref{1.15} instead of \eqref{1.14} and we proceed in the same way as above to show that the last inequality is still true by substituting $\|g-\tilde{g }\|_{L^2 (\Omega)}$ by $\|f-\tilde{f}\|_{H^1(\Omega)}$.

\section{Comments}\label{section3}

Throughout, the ball of $\mathbb{R}^n$ of radius $\rho >0$ and center $0$ will be denoted by $B_\rho$. Let $K$ be the closure of a bounded open subset of $\mathbb{R}^n$ with $C^\infty$ boundary. We say that  $K$  admits an escape function if there exists $B_\rho \supset K$ and $p\in C^\infty (\mathbb{S}^{n-1}\times \overline{U})$, where $U=B_\rho\setminus K$, satisfying : 
\\
(i) $\xi\cdot \partial_xp(\xi ,x)$ for each $(\xi,x)\in \mathbb{S}^{n-1}\times \overline{U}$.
\\
(ii) If $(\xi ,x)\in \mathbb{S}^{n-1}\times \partial K$ and $\eta =\xi -2(\xi\cdot \nu(x))\nu(x)$ then 
\[
p(\xi ,x)>p(\eta ,x)\; \mbox{when}\; \xi \cdot \nu (x)>0\quad \mbox{and}\quad \nu(x)\cdot \partial_\xi p(\xi ,x)>0\; \mbox{when}\; \xi \cdot \nu(x)=0.
\]
Here $\nu(x)$ denotes the exterior unit normal vector at $x$.

When $K$ is star-shaped with respect to $0$ we check that $p(\xi ,x)=\xi\cdot x$ is an escape function for $K$.

Below we assume that $K$ admits an escape function and we set  $\Omega =\mathbb{R}^n\setminus K$.

Fix $0<c_0\le 1<c_1$ and denote by $\mathcal{C}$ the set of functions $c\in C^\infty(\overline{\Omega})$ satisfying $c_0\le c\le c_1$ and $c=1$ outside $B_{\rho_0}\supset K$, for some fixed $\rho_0>0$.

Let $c\in \mathcal{C}$ and consider the unbounded operator $A_c:L^2(\Omega ,c^{-2}dx)\rightarrow L^2(\Omega ,c^{-2}dx)$ given by 
\[
A_cu=c^2\Delta u,\quad u\in D(A)=H_0^1(\Omega )\cap H^2(\Omega).
\]
We verify that $A_c$ enjoys from the same properties as the operator $A$ of Section \ref{section1}. Proceeding in the same way as in Section \ref{section1}, we show that for every $(f,g)\in (H_0^1(\Omega)\cap H^2(\Omega))\times H_0^1(\Omega)$ the IBVP 
\begin{equation}\label{3.1}
\partial_t^2u-c^2\Delta u=0\; \mbox{in}\; \mathbb{R}_+\times \Omega, \quad u_{|\mathbb{R}_+\times \partial \Omega}=0,\quad (u(0,\cdot),\partial_tu(0,t))=(f,g)
\end{equation}
admits a unique solution $u=u(c,f,g)$ so that $u\in C(\overline{\mathbb{R}}_+,H_0^1(\Omega)\cap H^2(\Omega))$, $\partial_tu\in C(\overline{\mathbb{R}}_+,H_0^1(\Omega))$ and $\partial_t^2u\in C(\overline{\mathbb{R}}_+,L^2(\Omega))$. Moreover, we have for every $t\ge 0$
\[
\|u(t,\cdot)\|_{H^2(\Omega)}+\|\partial_tu(t,\cdot)\|_{H^1(\Omega)}+\|\partial_t^2u(t,\cdot)\|_{L^2(\Omega)}\le \mathfrak{m} \left(\|f\|_{H^2(\Omega)}+\|g\|_{H^1(\Omega)}\right),
\]
where $\mathfrak{m}=\mathfrak{m}\left(\Omega, c_0,c_1\right)>0$ is a constant.

IBVPs (with $g=0$) of type \eqref{3.1} appear naturally in thermo-acoustic and photo-acoustic tomography, where $c$ represents the speed of sound and $f$ an internal source (see for example \cite{LU} and the references therein for more details).

According to \cite[Section 6]{MRS}, Theorem \ref{theorem1} still holds when \eqref{1.2} is replaced by~\eqref{3.1}. 

Let $c\in \mathcal{C}$ and $(f,g),(\tilde{f},\tilde{g})\in \mathbf{D}_\vartheta$, where $\vartheta >0$ is fixed. If $u=u(c,f,g)$ and $\tilde{u}=u(c,\tilde{f},\tilde{g})$. In this case, instead of \eqref{1.14} and \eqref{1.15} we have the following equations
\begin{align*}
&-\Delta w_0=c^{-2}(g-\tilde{g})\; \mbox{in}\;Q, \quad w_0{_{|\partial K}}=0, 
\\
&-\Delta w_1=c^{-2}(f-\tilde{f})\; \mbox{in}\;Q, \quad w_1{_{|\partial K}}=0.
\end{align*}
We can then mimic the proof of Theorem \ref{theorem1.0} in order to obtain when $c^{-2}(f-\tilde{f})$ and $c^{-2}(g-\tilde{g})$ are independent of $x_n$
\begin{equation}\label{mi2}
C\left(\|f-\tilde{f}\|_{L^2(\Omega)}+\|g-\tilde{g}\|_{L^2(\Omega)}\right)\le \rho^{-1}\vartheta +e^{e^{\mathfrak{p}\rho}}N(u-\tilde{u}),\quad \rho\ge 1,
\end{equation}
where $C=C(n,Q,S,c_0)>0$ and $\mathfrak{p}=\mathfrak{p}(n,Q,Q_0,c_0)>0$ are constants.

Next, let $c,\tilde{c}\in \mathcal{C}$ and $(f,g)\in \mathbf{D}_\vartheta$. Set $u=u(c,f,g)$ and $\tilde{u}=u(\tilde{c},f,g)$. In the same way as \eqref{mi2}, we obtain, assuming that $f(c^{-2}-\tilde{c}^{-2})$ and $g(c^{-2}- \tilde{ c}^{-2})$ are independent of $x_n$,\begin{equation}\label{mi3}
C\left(\|f(c^2-\tilde{c}^2)\|_{L^2(\Omega)}+\|g(c^2-\tilde{c}^2)\|_{L^2(\Omega)}\right)\le \rho^{-1}\vartheta +e^{e^{\mathfrak{p}\rho}}N(u-\tilde{u}),\quad \rho\ge 1,
\end{equation}
where $C$ and $\mathfrak{p}$ are as above.

Let us further assume that $f$ and $g$ are continuous. Let $\mathcal{K}_0$ be a compact subset of $\{f\ne 0\}$ and $\mathcal{K}_1$ be a compact subset of $\{g\ne 0\}$. If $\displaystyle m=\min \left(\min_{\mathcal{K}_0}|f|,\min_{\mathcal{K}_1}|g|\right)$ then \eqref{mi3} implies
\begin{equation}\label{mi4}
mC\|c-\tilde{c}\|_{L^2(\mathcal{K}_0\cup \mathcal{K}_1)}\le \rho^{-1}\vartheta +e^{e^{\mathfrak{p}\rho}}N(u-\tilde{u}),\quad \rho\ge 1,
\end{equation}
where $C$ and $\mathfrak{p}$ are the constants in \eqref{mi3}.

As it was already observed in \cite{KU,LU},  if one wants to  simultaneously recover the initial data $f$ (in the case $g=0$) and the  speed of sound $c$ from a single measurement then the best we can do is to recover $c^{-2}f$ under an additional hypothesis. 

A series of results were obtained by Stefanov and Uhlmann \cite{SU1,SU2,SU3} for the problem of determining $c$ from the map $(f,0)\mapsto u(c,f,0)_{|(0,T)\times \Sigma}$ with $T\gg 1$. The method employed in \cite{SU1,SU2,SU3} is based on the direct linearization of the wave equation. Under different conditions used to obtain \eqref{mi4}, the authors prove in \cite{SU3} a Lipschitz stability inequality when $\mbox{supp}(c-\tilde{c})$ is contained in $\{\Delta f\ne 0\}$, $f$ sufficiently smooth and $g=0$.

In \cite{SU4} Stefanov and Uhlmann show that the inverse problem mentioned above is unstable at any scale of Sobolev spaces.

The relationship between this type of inverse problems and the transmission eigenvalues was established in \cite{FH} by Fink and Hickmann.

Let us now define
\[
\mathscr{D}=\{(f,g)\in (H_0^1(\Omega) \cap H^2(\Omega))\cap H_0^1(\Omega);\; f\; \mbox{and}\; g\; \mbox{have compact support}\}.
\]
We consider, where $c\in \mathcal{C}$, the following two initial-to-boundary operators
\begin{align*}
&\Phi_c^0:(f,0)\in \mathscr{D}\mapsto (u(c,f,0),\nabla u(c,f,0)_{|\mathbb{R}_+\times S}\in L^1(\mathbb{R}_+,L^2(S))^{n+1},
\\
&\Phi_c^1:(0,g)\in \mathscr{D}\mapsto (u(c,0,g),\nabla u(c,0,g)_{|\mathbb{R}_+\times S}\in L^1(\mathbb{R}_+,L^2(S))^{n+1}.
\end{align*}

The following corollary is a direct consequence of inequality \eqref{mi3}, where $\mathcal{C}_j=\{c\in \mathcal{C};\; c\; \mbox{is independent of  the variable}\; x_j\}$, $1\le j\le n$,

\begin{corollary}\label{corollary3.1}
$\Phi_c^0$ or $\Phi_c^1$ determines uniquely the speed of sound $c\in \mathcal{C}_j$, $1\le j\le n$.
\end{corollary}

\end{document}